\documentclass[12pt]{amsart}
\usepackage{amssymb}
\usepackage{fullpage}
\usepackage{lscape}
\usepackage{graphicx}

\newcommand{\frp}{\mathfrak p}
\newcommand{\frq}{\mathfrak q}

\newcommand{\ddim}{\mathbf{dim}}

\newcommand{\sub}{\subseteq}

\newcommand{\Z}{\mathbb Z}

\renewcommand{\a}{\mathbf a}
\renewcommand{\d}{\mathbf d}

\newcommand{\cA}{\mathcal A}
\newcommand{\cF}{\mathcal F}
\newcommand{\cM}{\mathcal M}
\newcommand{\cO}{\mathcal O}
\newcommand{\cQ}{\mathcal Q}

\newcommand{\rev}{\mathrm{rev}}

\DeclareMathOperator{\GL}{GL} %

\newtheorem{thm}{Theorem}[section]
\newtheorem{lem}[thm]{Lemma}

\theoremstyle{remark}
\newtheorem{rem}[thm]{Remark}
\newtheorem{rems}[thm]{Remarks}

\numberwithin{equation}{section}

\title[Orbits of parabolic subgroups on metabelian ideals]
{Orbits of parabolic subgroups on metabelian ideals}

\author[S.~M.~Goodwin]{Simon M.~Goodwin}
\address{School of Mathematics, University of Birmingham, Birmingham, B15
2TT, UK}  \email{goodwin@maths.bham.ac.uk}
\urladdr{http://web.mat.bham.ac.uk/S.M.Goodwin/}

\author[L.~Hille]{Lutz Hille}
\address{Fakult\"at f\"ur Mathematik, Universit\"at Bielefeld, Postfach 100131, 33501 Bielefeld, Germany}
\email{hille@math.uni-bielefeld.de}
\urladdr{http://www.math.uni-bielefeld.de/$\sim$hille/}

\author[G.~R\"ohrle]{Gerhard R\"ohrle}
\address{Fakult\"at f\"ur Mathematik, Ruhr-Universit\"at Bochum,
Universit\"atsstrasse 150, D-44780 Bochum, Germany}
\email{gerhard.roehrle@rub.de}
\urladdr{http://www.ruhr-uni-bochum.de/ffm/Lehrstuehle/Lehrstuhl-VI/rubroehrle.html}

\thanks{2000 {\it Mathematics Subject Classification}: 14L30, 20G15}

\begin{document}

\begin{abstract}
Let $k$ be an algebraically closed field, $t \in \Z_{\ge1}$, and let $B$ be the Borel
subgroup of $\GL_t(k)$ consisting of upper-triangular matrices. Let
$Q$ be a parabolic subgroup of $\GL_t(k)$ that contains $B$ and such
that the Lie algebra $\frq_u$ of the unipotent radical of $Q$ is
metabelian, i.e.\ the derived subalgebra of $\frq_u$ is abelian. For
a dimension vector $\d = (d_1,\dots,d_t) \in \Z_{\ge 1}^t$ with
$\sum_{i=1}^t d_i = n$, we obtain a parabolic subgroup $P(\d)$ of
$\GL_n(k)$ from $B$ by taking upper-triangular block matrices with
$(i,j)$ block of size $d_i \times d_j$.  In a similar manner we
obtain a parabolic subgroup $Q(\d)$ of $\GL_n(k)$ from $Q$.  We
determine all instances when $P(\d)$ acts on $\frq_u(\d)$ with a
finite number of orbits for all dimension vectors $\d$. Our methods
use a translation of the problem into the representation theory of
certain quasi-hereditary algebras. In the finite cases, we use
Auslander--Reiten theory to explicitly determine the $P(\d)$-orbits;
this also allows us to determine the degenerations of
$P(\d)$-orbits.
\end{abstract}

\maketitle

\section{Introduction}

Let $G$ be a reductive algebraic group over an algebraically closed
field $k$, and assume that the characteristic of $k$ is good for
$G$. There has been a great deal of recent interest in the adjoint
action of a parabolic subgroup $P$ of $G$ on the Lie algebra
$\frp_u$ of its unipotent radical. Such actions were considered in
the case where $\frp_u$ is abelian by R.~Richardson, R.~Steinberg
and the third author in \cite{RRS}; in particular, they showed that
if $\frp_u$ is abelian, then there is always a finite number of
$P$-orbits in $\frp_u$ and gave a parameterization of these orbits.
Subsequently, in work of U.~J\"urgens and the second and third
authors, a classification of all instances when $P$ acts on $\frp_u$
with finitely many orbits has been obtained, see \cite{HR} and
\cite{JR}. In addition, for $G$ simple not of type $E_7$ or $E_8$
there is a classification of all instances when there are finitely
many $P$-orbits in higher terms $\frp_u^{(l)}$ of the descending
central series of $\frp_u$, see \cite{BH}, \cite{BHR} and \cite{GR}.
For further information on parabolic group actions we refer the
reader to the survey \cite{Romono}.

In the case $G = \GL_n(k)$, there has been much success in
understanding the adjoint action of a parabolic subgroup through a
translation in to the representation theory of certain
quasi-hereditary algebras. This translation was first observed in
\cite{HR}, and has subsequently been further exploited, see for
example, \cite{BH} and \cite{BHRZ}.  We refer the reader to \cite{H}
and \cite{JSY} for recent related developments.

In this paper we consider a related problem in case $G = \GL_n(k)$.
Rather than considering the action of $P$ on $\frp_u$ (or
$\frp_u^{(l)}$), we study the action of $P$ on $\frq_u$, where
$Q$ is a parabolic subgroup of $G$ such that $\frq_u$ is metabelian.
Our main result, Theorem \ref{T:main}, gives a finiteness condition
for such actions.

\smallskip

Let $G = \GL_n(k)$, where $k$ is an algebraically closed field. Let
$\d = (d_1,\dots,d_t) \in \Z_{\ge 1}^t$ satisfy $\sum_{i=1}^t d_i =
n$; call such a $t$-tuple $\d$ a {\em dimension vector}. Given a dimension
vector $\d$, define the parabolic subgroup $P = P(\d)$ of $G$ to
be the stabilizer of the {\em standard flag} $0 \sub k^{e_1} \sub
k^{e_2} \sub \dots \sub k^{e_{t-1}} \sub k^{e_t} = k^n$ in $k^n$, where $e_i
= \sum_{j=1}^i d_j$. Let $\a = (a_1,a_2,a_3) \in \Z_{\ge 1}^3$ with
$a_1+a_2+a_3 =t$ and define $b_i = \sum_{j=1}^i a_j$ for $i =
0,1,2,3$. Define $Q = Q(\a,\d) = P(e_{b_1}, e_{b_2}-e_{b_1},
e_{b_3} - e_{b_2}) \supseteq P$. Then $P$ acts on the Lie algebra
$\frq_u = \frq_u(\a,\d)$ of the unipotent radical of $Q$ via the
adjoint action.

Define the partial order $\le$ on the set of triples $\a \in \Z_{\ge
1}^3$, by $\a \le \tilde \a$ if $a_i \le \tilde a_i$, for $i =
1,2,3$.  The {\em reverse} of a triple $\a$ is 
$\a_\rev = (a_3,a_2,a_1)$.

We can now state our main result which gives 
a classification of all triples
$\a$ such that $P(\d)$ acts with a finite number of orbits on
$\frq_u(\a,\d)$ for all dimension vectors $\d$.  


\begin{thm} 
\label{T:main}
Let $\a \in \Z_{\ge 1}^3$ and $\d \in \Z_{\ge 1}^t$ a dimension vector  
with $a_1+a_2+a_3 =t$.
\begin{enumerate}
\item Assume $\a$ and $\d$ are as in Table \ref{Tab:mininf}.
Then $P(\d)$ acts on $\frq_u(\a,\d)$ with an infinite number of
orbits.
\item Assume $\a$ is a triple in Table
\ref{Tab:maxfin}.  Then $P(\d)$ acts on $\frq_u(\a,\d)$ with a
finite number of orbits for all $\d$.
\item Any triple $\a$ is either: less than or equal to a triple in Table
\ref{Tab:maxfin} or its reverse; or greater than or equal to a
triple in Table \ref{Tab:mininf} or its reverse. In the former case
there is a finite number of $P(\d)$-orbits in $\frq_u(\a,\d)$ for
all $\d$, and in the latter case there is an infinite number of
$P(\d)$-orbits in $\frq_u(\a,\d)$ for some $\d$.
\end{enumerate}
\end{thm}

\renewcommand{\arraystretch}{1.2}
\begin{table}[hbt]
\begin{tabular}{|c|c|c|}
\hline $\a$ & $\d$ & $\dim P/Q_u$ \\
\hline \hline $(1,3,5)$ & $(3,2,2,2,1,1,1,1,1)$ & $48$\\
\hline $(1,4,3)$ & $(2,1,1,1,1,1,1,1)$ & $20$\\
\hline $(1,6,2)$ & $(3,1,1,1,1,1,1,2,2)$ & $42$\\
\hline $(2,2,5)$ & $(2,2,3,3,1,1,1,1,1)$ & $54$\\
\hline $(2,3,2)$ & $(1,1,1,1,1,1,1)$ & $12$\\
\hline $(3,2,3)$ & $(1,1,1,2,2,1,1,1)$ & $24$\\
\hline
\end{tabular}
\bigskip \caption{The minimal infinite triples $\a$}
\label{Tab:mininf}
\end{table}
\begin{table}[hbt]
\begin{tabular}{|c||c|c|c|c|c|c|}
\hline $\a$ & $(1,2,a)$ & $(1,a,1)$ &  $(a,1,b)$ & $(1,3,4)$ &
$(1,5,2)$ & $(2,2,4)$
\\ \hline
\end{tabular}
\bigskip \caption{The maximal finite triples $\a$ ($a,b \in \Z_{\ge1}$)} \label{Tab:maxfin}
\end{table}

\begin{rems}
(i).
Observe that $P(\d)$ has a finite number of orbits on
$\frq_u(\a,\d)$ for all $\d$ if and only if the same is true for the
action of $P(\d)$ on $\frq_u(\a_\rev,\d)$. This easy observation
means that Theorem \ref{T:main} does indeed give the desired
classification of all triples $\a$ such that $P(\d)$ acts on
$\frq_u(\a,\d)$ with a finite number of orbits for all $\d$.

(ii).
Note that Theorem \ref{T:main} is consistent with the
computer calculations made in \cite{HV}.

(iii).
Compared with Theorem \ref{T:main}, 
it is a much harder problem to determine all pairs
$(\a,\d)$ such that $P(\d)$ acts on $\frq_u(\a,\d)$ with a finite
number of orbits.  For the minimal infinite cases given in Table
\ref{Tab:mininf}, we have that $P(\d)$ acts on $\frq_u(\a,\d)$ with
finitely many orbits if one entry of $\d$ is less than that for the
$\d$ given in the table.  For larger values of $\a$ it can also be
the case that $P(\d)$ acts on $\frq_u(\a,\d)$ with finitely many
orbits for some values of $\d$.  It seems infeasible to determine
all such $\d$ for all $\a \in \Z_{\ge 1}^3$.
\end{rems}

\smallskip

We prove Theorem \ref{T:main} in the next section.  The proof of (1)
is a dimension counting argument that is elementary; it essentially
involves calculating the dimension given in the third column of
Table \ref{Tab:mininf}. The majority of the work required is in
proving (2). To do this we interpret the problem in terms of the
representation theory of a certain quasi-hereditary algebra
$\cA(\a)$.  The isoclasses of $\Delta$-filtered $\cA(\a)$-modules
with $\Delta$-dimension vector $\d$ correspond to the orbits of
$P(\d)$ in $\frq_u(\a,\d)$.  For each of the triples in Table
\ref{Tab:maxfin}, we have calculated the Auslander--Reiten quiver of
$\Delta$-filtered $\cA(\a)$-modules. In doing so we see that for
these values of $\a$ the algebra $\cA(\a)$ has finite
$\Delta$-representation type, which proves (2).

In order to prove (3) we just require the following observation. If
$\tilde \a \le \a$, then there is an embedding of the category of
$\Delta$-filtered $\cA(\tilde \a)$-modules in to the category of
$\Delta$-filtered $\cA(\a)$-modules.  The correspondence alluded to
above then implies that for $\tilde \a \le \a$, if there are
finitely $P(\d)$-orbits in $\frq_u(\a,\d)$ for all $\d$, then there
are finitely many $P(\tilde \d)$-orbits in $\frq_u(\tilde \a,\tilde
\d)$ for all $\tilde \d$.

\begin{rem}
Using the results from \cite{BHRZ} and the Auslander--Reiten quivers
that we have calculated, it is possible to calculate the
degenerations of $P(\d)$-orbits in $\frq_u(\a,\d)$ in the finite
cases.  More precisely, let $\cO$, $\cO'$ be $P$-orbits and $M$,
$M'$ the corresponding $\cA(\a)$-modules.  In \cite{BHRZ} it is
shown that $M$ is greater than $M'$ in the {\em hom-order} if and
only if $\cO$ is a degeneration of $\cO'$. The hom-order can be read
off from the Auslander--Reiten quiver.
\end{rem}



\section{Proof of Theorem \ref{T:main}} \label{S:proof}

Let $\a \in \Z_{\ge 1}^3$ and $\d \in \Z_{\ge 1}^t$ with $a_1+a_2+a_3 =t$.  
Let $P = P(\d)$ and $Q = Q(\a,\d)$.  Write $Q_u = Q_u(\a,\d)$ for the
unipotent radical of $Q$ and $\frq_u'= \frq_u'(\a,\d)$ for the
derived subalgebra of $\frq_u= \frq_u(\a,\d)$.  Define $b_i =
\sum_{j=1}^i a_j$ for $i = 0,1,2,3$ as in the introduction.

\subsection*{Proof of (1)}
The action of $P$ on $\frq_u$ induces an action of $P/Q_u$ on
$\frq_u/\frq_u'$, and if $P$ acts on $\frq_u$ with finitely many
orbits, then there is necessarily a finite number of $P/Q_u$-orbits
in $\frq_u/\frq_u'$.  For the dimension vectors given in the second
column of Table \ref{Tab:mininf}, we have $\dim P/Q_u = \dim
\frq_u/\frq_u'$; this dimension is given in the third column.
Observe that the scalar matrices in $P/Q_u$ act trivially on
$\frq_u/\frq_u'$.  Therefore, there cannot be a dense $P/Q_u$-orbit
on $\frq_u/\frq_u'$ and hence there must be infinitely many
$P/Q_u$-orbits.  This proves part (1) of Theorem \ref{T:main}.

\begin{rem}
The above proof of (1) is obtained by considering a particluar quadratic
form in the dimension vectors $\d \in \Z_{\ge 1}^t$. 
For a fixed $\a = (a_1, a_2, a_3)$ with $a_1+a_2+a_3 =t$ and $b_i$ ($0 \le i \le 3$) as before, 
set $I = \{(i,j) \mid b_{k-1} < i < j \le b_k \text{ for some }
k=1,2,3\}$ and $J = \{(i,j) \mid b_{k-1} < i \le b_k < j \le b_{k+1}
\text{ for some } k=1,2\}$. 
The quadratic form at $\d \in \Z_{\ge 1}^t$ is then given by  
$$
\sum_{i=1}^t d_i^2 + \sum_{(i,j) \in I} d_i d_j - \sum_{(i,j) \in J}
d_i d_j.
$$
On readily checks that, by our construction, this expression 
equals $\dim P/Q_u - \dim \frq_u/\frq_u'$.

For each of the triples $\a$ in Table \ref{Tab:mininf}, this
quadratic form is semi-definite and the value of $\d$ given in the
table is the unique indivisible vector in the kernel of the
quadratic form. For a classification of critical quadratic unit
forms we refer to \cite{HHoehne} and for further applications in
representation theory to Ringel's monograph \cite{RiTame}.
\end{rem}

\subsection*{Proof of (2)}
In order to prove part (2) of Theorem \ref{T:main}, we translate the
problem in to one regarding the representation theory of a certain
quasi-hereditary algebra $\cA(\a)$.  The algebra $\cA(\a)$ is
defined as the quotient of a path algebra by some relations as
explained below.

The quiver $\cQ = \cQ(\a)$ is defined to have vertex set
$\{1,\dots,t\}$; there are arrows $\alpha_i : i \to i+1$ for all
$i$, an arrow $\beta_{b_1} : b_2 \to b_1$ and an arrow $\beta_{b_2}
: b_3 \to b_2$. Below, in Figure \ref{F:Q}, we give an example of a
quiver $\cQ(\a)$.
\begin{figure}[h!]
\begin{picture}(320,25)
\put(0,15){\circle*{5}} %
\put(0,5){\makebox(0,0){\tiny 1}} %
\put(40,15){\circle*{5}} %
\put(40,5){\makebox(0,0){\tiny 2}} %
\put(80,15){\circle*{5}} %
\put(80,5){\makebox(0,0){\tiny 3}} %
\put(120,15){\circle*{5}} %
\put(120,5){\makebox(0,0){\tiny 4}} %
\put(160,15){\circle*{5}} %
\put(160,5){\makebox(0,0){\tiny 5}} %
\put(200,15){\circle*{5}} %
\put(200,5){\makebox(0,0){\tiny 6}} %
\put(240,15){\circle*{5}} %
\put(240,5){\makebox(0,0){\tiny 7}} %
\put(280,15){\circle*{5}} %
\put(280,5){\makebox(0,0){\tiny 8}} %
\put(320,15){\circle*{5}} %
\put(320,5){\makebox(0,0){\tiny 9}} %
\put(5,10){\vector(1,0){30}} %
\put(20,0){\makebox(0,0){\small $\alpha_1$}} %
\put(45,10){\vector(1,0){30}} %
\put(60,0){\makebox(0,0){\small $\alpha_2$}} %
\put(85,10){\vector(1,0){30}} %
\put(100,0){\makebox(0,0){\small $\alpha_3$}} %
\put(125,10){\vector(1,0){30}} %
\put(140,0){\makebox(0,0){\small $\alpha_4$}} %
\put(165,10){\vector(1,0){30}} %
\put(180,0){\makebox(0,0){\small $\alpha_5$}} %
\put(205,10){\vector(1,0){30}} %
\put(220,0){\makebox(0,0){\small $\alpha_6$}} %
\put(245,10){\vector(1,0){30}} %
\put(260,0){\makebox(0,0){\small $\alpha_7$}} %
\put(285,10){\vector(1,0){30}} %
\put(300,0){\makebox(0,0){\small $\alpha_8$}} %
\put(155,20){\vector(-1,0){70}} %
\put(120,30){\makebox(0,0){\small $\beta_3$}} %
\put(315,20){\vector(-1,0){150}} %
\put(240,30){\makebox(0,0){\small $\beta_5$}} %
\end{picture}
\caption{The quiver $\cQ(3,2,4)$} \label{F:Q}
\end{figure}

\noindent Let $I$ be the ideal of the path algebra $k\cQ$ of
$\cQ$ generated by the relations:
\begin{equation} 
\label{e:relns}
\beta_{b_1} \alpha_{b_2-1} \cdots \alpha_{b_1+1} \alpha_{b_1} = 0
\quad \text{ and } \quad \alpha_{b_2-1} \cdots \alpha_{b_1+1} \alpha_{b_1}
\beta_{b_1} = \beta_{b_2} \alpha_{b_3-1} \cdots \alpha_{b_2+1}
\alpha_{b_2}.
\end{equation}
The algebra $\cA = \cA(\a)$ is defined to be the quotient $k\cQ/I$.

Recall that an $\cA$-module is given by a family of vector spaces
$M_i$ for $i = 1,\dots,t$ and linear maps $M_{\alpha_i} : M_i \to
M_{i+1}$ for $i = 1,\dots,t-1$, and $M_{\beta_{b_i}} : M_{b_{i+1}}
\to M_{b_i}$ for $i = 1,2$ that satisfy the relations
\eqref{e:relns}. For an arrow $\gamma$ in $\cQ$ we often simply write
$\gamma$ rather than $M_\gamma$ when considering an $\cA$-module
$M$.

Note that the arguments of \cite[\S 2--3]{BH} apply in the
present situation showing that $\cA$ is a quasi-hereditary algebra,
and that the isoclasses of $\Delta$-filtered $\cA$-modules with
$\Delta$-dimension vector $\d$ are in bijective correspondence with
the adjoint orbits of $P = P(\d)$ on $\frq_u = \frq_u(\a,\d)$.

Define the
category $\cM(\a)$ to be the full subcategory of $\cA(\a)$-mod
consisting of modules $M$ for which $M_{\alpha_i}$ is injective for
all $i = 1,\dots,t-1$.  Let $\cM(\a,\d)$ be the class of modules
$M$ in $\cM(\a)$ for which $\dim M_i = e_i = \sum_{j=1}^i d_j$ for
all $i = 1,\dots,t$.
Our next result can be proved in exactly
the same way as \cite[Lem.\ 1 and Lem.\ 2]{BH}.  

\begin{lem} 
\label{L:equiv}
The adjoint orbits of $P(\d)$ on $\frq_u(\a,\d)$ are in bijective
correspondence with the isoclasses of modules in $\cM(\a,\d)$.
\end{lem}

Define the {\em standard $\cA(\a)$-module} $\Delta(i)$ as follows. For $j <
i$, we have $\Delta(i)_j = 0$ and for $j \ge i$, we have
$\Delta(i)_j = k$. The arrows $\alpha_j$ act as the zero map for $j < i$
and as the identity map for $j \ge i$; the arrows $\beta_{b_1}$
and $\beta_{b_2}$ both act as zero maps. With these definitions one can
prove the following lemma in exactly the same way as \cite[Prop.\
1]{BH}.

\begin{lem} \label{L:qh}
The algebra $\cA(\a)$ is quasi-hereditary with the inverse
order on the index set $\{1,\dots,t\}$ and standard modules as
defined above. The category $\cM(\a)$ is precisely the category
$\cF(\cA(\a),\Delta)$ of $\Delta$-filtered $\cA(\a)$-modules.
\end{lem}

Recall that an $\cA$-module $M$ is called {\em $\Delta$-filtered}
if there exists a filtration $0 = M_0 \sub M_1 \sub \dots \sub M_r =
M$ of $M$, where for each $i$ we have $M_i/M_{i-1} \cong \Delta(j)$
for some $j$; such a filtration of $M$ is called a {\em
$\Delta$-filtration}. Given a $\Delta$-filtered module $M$, define
the {\em $\Delta$-dimension vector} $\ddim_\Delta M \in \Z_{\ge
0}^t$ of $M$ by setting $(\ddim_\Delta M)_j$ equal to the number of
factors isomorphic to $\Delta(j)$ in a $\Delta$-filtration of $M$;
it is a standard result that the $\Delta$-dimension vector is
well-defined. Recall that a quasi-hereditary algebra is said to
have {\em finite $\Delta$-representation type} if there are only
finitely many isoclasses of indecomposable $\Delta$-filtered
modules. For general results on categories of $\Delta$-filtered
modules over quasi-hereditary algebras, we refer to \cite{DR}.

It follows from Lemmas \ref{L:equiv} and \ref{L:qh} that the
$P(\d)$-orbits in $\frq_u(\a,\d)$ are in bijective correspondence with
the $\Delta$-filtered $\cA(\a)$-modules with $\Delta$-dimension
vector $\d$.  Using this equivalence we verify part (2) of Theorem
\ref{T:main} by showing that $\cA(\a)$ has finite
$\Delta$-representation type for each triple $\a$ in Table
\ref{Tab:maxfin}.  We achieve this by calculating the
Auslander--Reiten quiver of $\cF(\cA(\a),\Delta)$, for each such
$\a$. For the values of $\a$ involving a parameter, we calculate the
Auslander--Reiten quiver for a particular value of that parameter;
it is then straightforward to generalize to the generic case.
For completeness, 
the Auslander--Reiten quivers are given in the appendix;
next we give a sketch of their construction. 

Thanks to \cite{Ri}, the category of
$\Delta$-filtered modules for a quasi-hereditary algebra has almost
split sequences. Therefore, the 
Auslander--Reiten quiver is defined for $\cF(\cA(\a),\Delta)$. 
To calculate these quivers, we use standard methods as explained in
\cite[Ch.\ VII]{ARS}, see also \cite{RiTame}. We begin by taking a
$\Z$-covering of the algebra $\cA(\a)$ and then making a ``cut''; we
denote the resulting algebra by $\cA_\Z(\a)$.  The category
$\cF(\cA_\Z(\a),\Delta)$ admits a simple projective object.
Thus we can construct the AR-quiver for
$\cF(\cA_\Z(\a),\Delta)$ by ``knitting''. After several steps 
repetitions occur and as a result, one can obtain the entire AR-quiver
for $\cF(\cA(\a),\Delta)$.

\subsection*{Proof of (3)}  To verify
the first assertion of (3), it suffices to check that if we obtain
the triple $\tilde \a$ from a triple $\a$ in Table \ref{Tab:mininf}
by subtracting 1 from an entry of $\a$, then $\tilde \a$ is less
than or equal to one of the triples in Table \ref{Tab:maxfin}.  This
is an elementary case by case check and is left to the reader.

\smallskip

Let $\tilde \a \le \a$.  We show below that there is embedding of
$\cF(\cA(\tilde \a),\Delta)$ into $\cF(\cA(\a),\Delta)$.  Along with
the correspondence between $P(\d)$-orbits in $\frq_u(\a,\d)$ and
$\Delta$-filtered $\cA(\a)$-modules with $\Delta$-dimension vector
$\d$ this implies the second assertion of (3).  The embedding
identifies $\cF(\cA(\tilde \a),\Delta)$ with the subcategory of
$\cF(\cA(\a),\Delta)$ consisting of modules with $\Delta$-dimension
vector $\d$ such that $d_i = 0$ for $i = \tilde a_1+1,\dots,b_1$, $i
= b_1+\tilde a_2+1,\dots,b_2$ and $i = b_2+\tilde a_3+1,\dots,b_3$.

In order to define the embedding, we require some notation. Define
$\tilde b_i$ for $i=0,1,2,3$ in analogy to the definition of $b_i$.
The map $f : \{1,\dots,\tilde b_3\} \to \{1,\dots,b_3\}$ is defined
by setting:
$$
f(i) = \left\{ \begin{array}{llrcl}
i - b_j + \tilde b_j & \text{if } & b_j  & < i \le & b_j + \tilde a_{j+1}; \\
\tilde b_{j+1} & \text{if} & b_j +\tilde a_{j+1} & < i \le & b_{j+1}.
\end{array} \right.
$$

Let $\tilde M$ be a $\Delta$-filtered $\cA(\tilde \a)$-module.  We
obtain a $\Delta$-filtered $\cA(\a)$-module $M$ from $\tilde M$ as
follows. Set $M_i = \tilde M_{f(i)}$ for $i = 1,\dots,t$. If
$f(i+1) = f(i)$, then define $M_{\alpha_i}$ to be the identity
map; and if $f(i+1) = f(i) + 1$, then define $M_{\alpha_i} =
\tilde M_{\alpha_{f(i)}}$. Finally, define, $M_{\beta_{b_i}} =
\tilde M_{\beta_{\tilde b_i}}$ for $i = 1,2$. This construction
defines a functor from $\cF(\cA(\tilde \a),\Delta)$ to
$\cF(\cA(\a),\Delta)$ by $\tilde M \mapsto M$.  It is
straightforward to check that this an embedding. This completes
the proof of (3).

\appendix

\section{Auslander--Reiten quivers}

\begin{figure}[h!b]
\centering
\includegraphics[width=14cm]{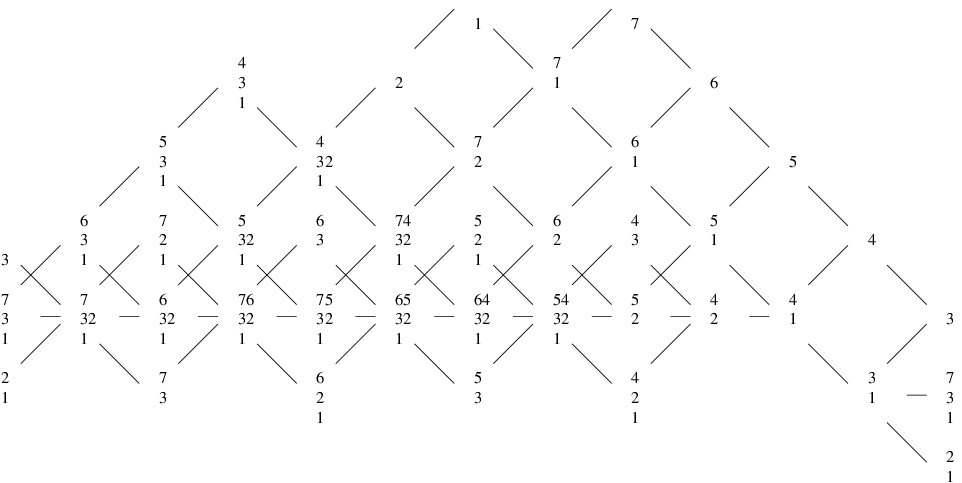}
\caption{The AR-quiver of $\Delta$-filtered modules
for $\cA(1,2,a)$ for $a=4$}  \label{Fig:ARfirst}
\end{figure}

\begin{landscape}
\begin{figure}
\centering
\includegraphics[width=22cm]{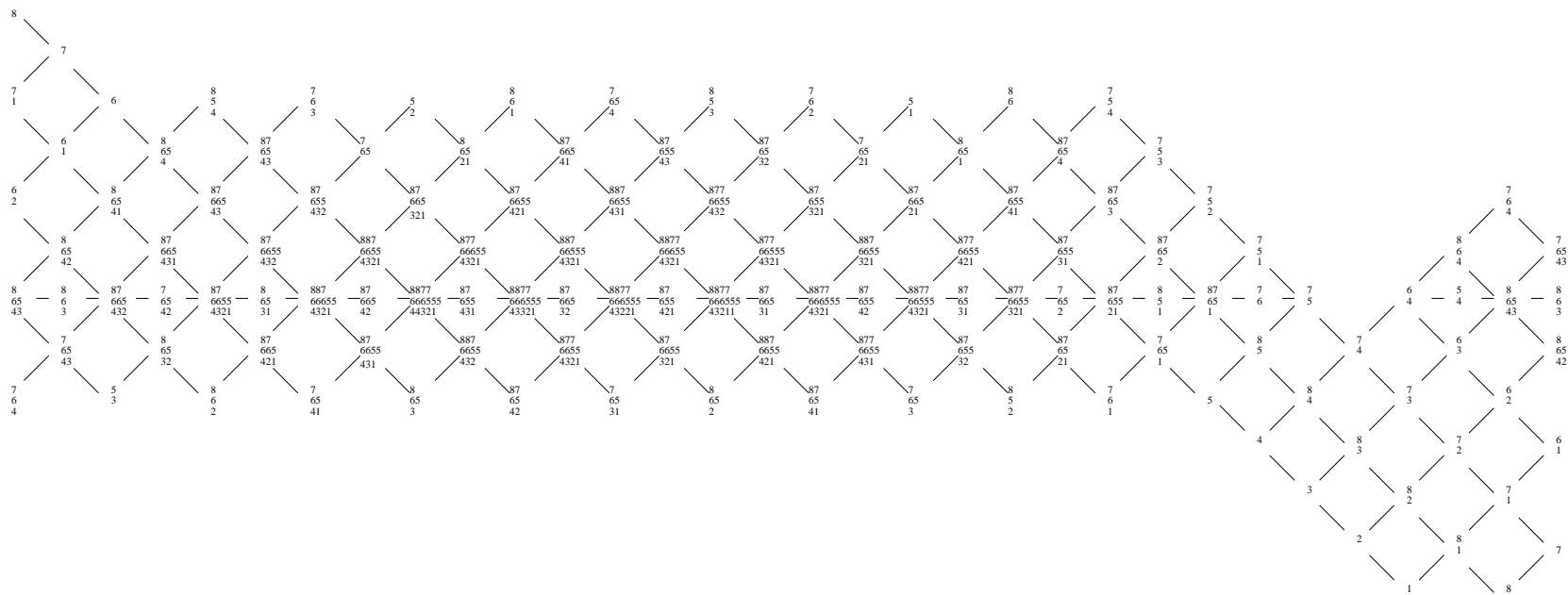}
\caption{The Auslander--Reiten quiver of $\Delta$-filtered modules
for $\cA(4,2,2)$} \label{Fig:ARlast}
\end{figure}
\end{landscape}

In Figures \ref{Fig:ARfirst}--\ref{Fig:ARlast} we include the
Auslander--Reiten quivers of the $\Delta$-filtered modules for the
algebras $\cA(\a)$ in case of finite $\Delta$-representation type.
The $\cA(\a)$-modules are denoted by giving their filtration into
standard modules. The right and left boundary of the
Auslander--Reiten quivers must be identified, it is a cylinder or a
M\"obius band.

\begin{figure}
\centering
\includegraphics[width=13cm]{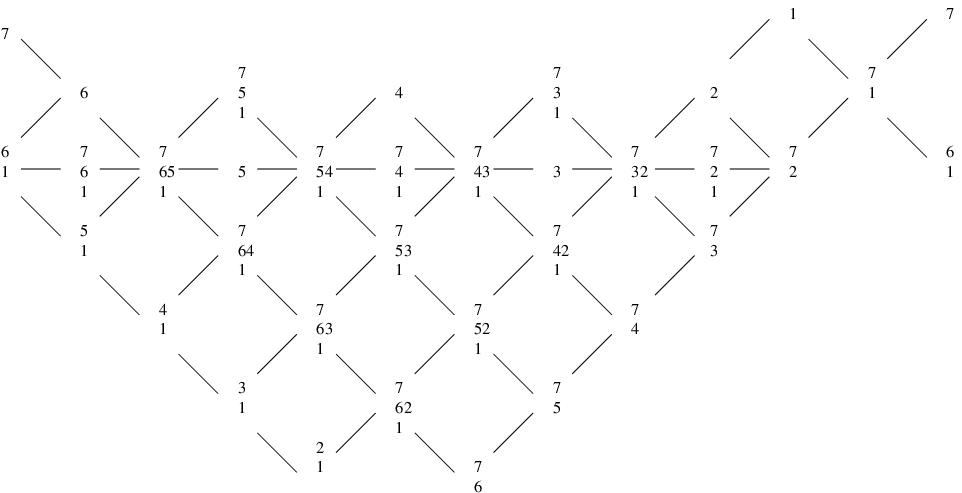}
\caption{The AR-quiver of $\Delta$-filtered modules
for $\cA(1,a,1)$ for $a = 5$} \vspace{20mm}
\end{figure}

\begin{figure}
\centering
\includegraphics[height=11cm]{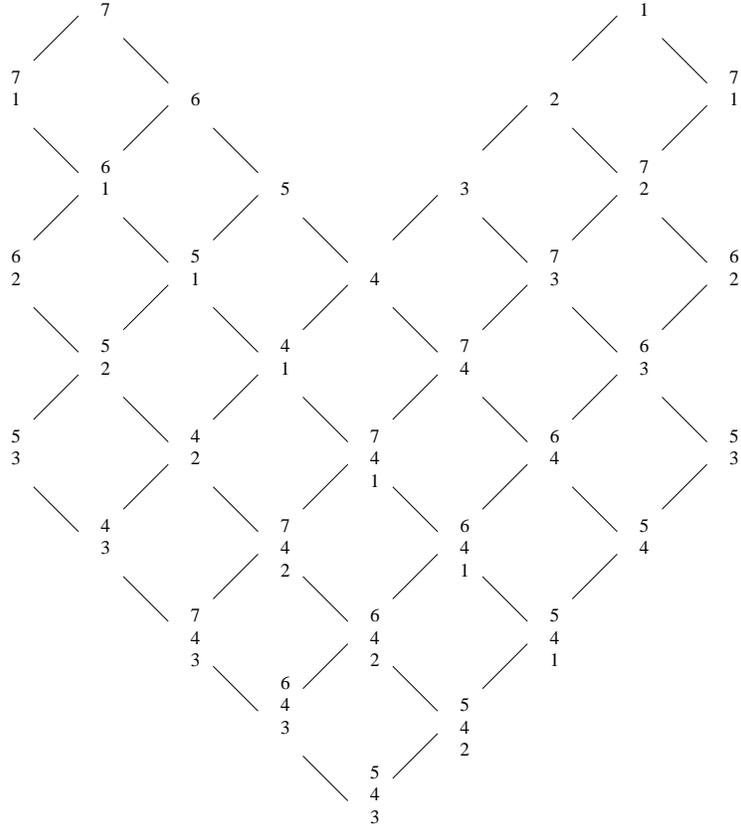}
\caption{The AR-quiver of $\Delta$-filtered modules
for $\cA(a,1,b)$ for $a = 3$, $b = 3$}
\end{figure}

\begin{landscape}
\begin{figure}
\centering
\includegraphics[width=22cm]{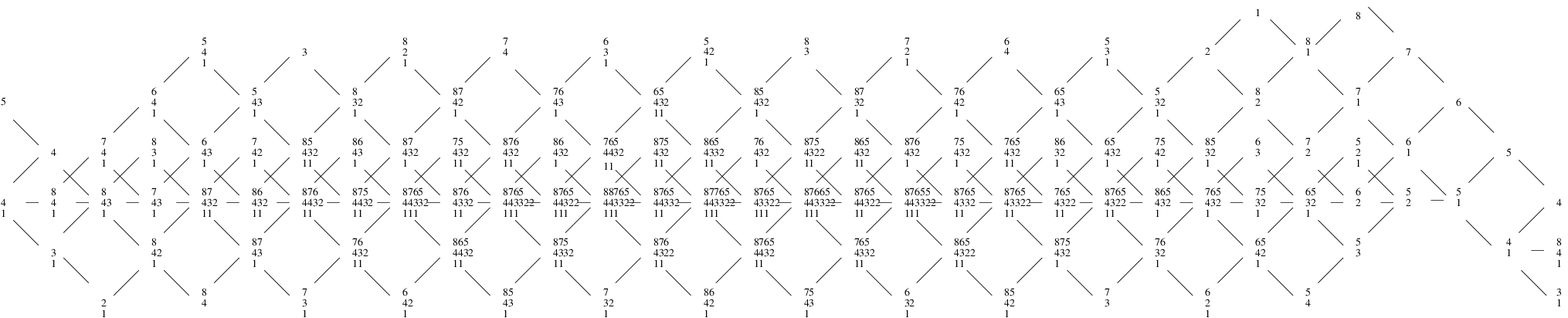}
\caption{The Auslander--Reiten quiver of $\Delta$-filtered modules
for $\cA(1,3,4)$} \vspace{20mm}
\end{figure}

\begin{figure}
\centering
\includegraphics[width=22cm]{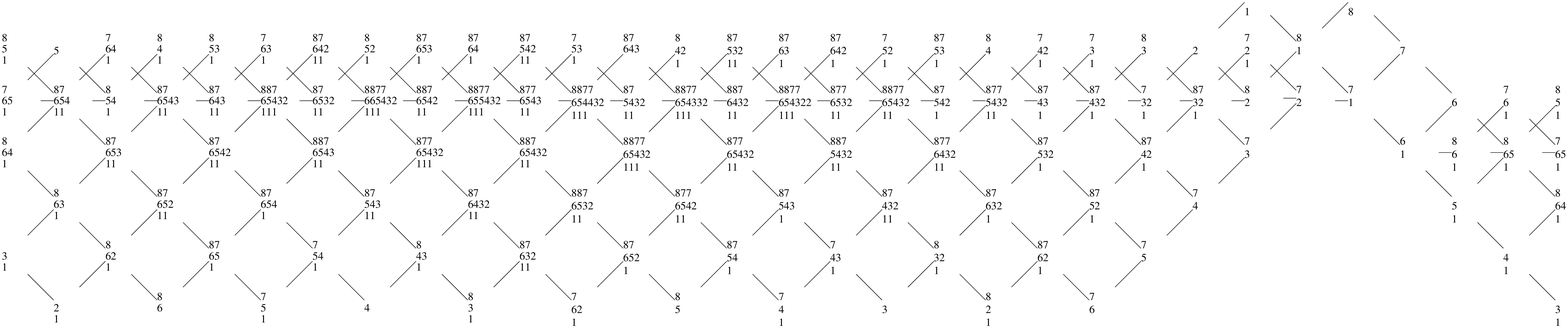}
\caption{The Auslander--Reiten quiver of $\Delta$-filtered modules
for $\cA(1,5,2)$}
\end{figure}
\end{landscape}

\section*{Acknowledgements}
This research was funded in part by EPSRC grant EP/D502381/1. Part
of the research for this paper was carried out while the authors
were staying at the Mathematical Research Institute Oberwolfach
supported by the "Research in Pairs" programme.

\end{document}